\documentclass[10pt]{amsart}
\usepackage{latexsym}
\usepackage{amssymb}
\usepackage{amsmath}
\usepackage{mathrsfs}
\usepackage{bbm}
\usepackage{bbold}
\def\Aff{\operatorname{Aff}}

\def\Fd{\operatorname{Fd}}
\def\spt{\operatorname{supp}}

\usepackage[bookmarks=true]{hyperref}
\usepackage{graphicx}

\begin{document}
\title{Dido's Problem and Beyond}

\author{S.~S. Kutateladze}
\address[]{
Sobolev Institute of Mathematics\newline
\indent 4 Koptyug Avenue\newline
\indent Novosibirsk, 630090\newline
\indent Russia}
\email{
sskut@math.nsc.ru
}
\begin{abstract}
We overview the main ideas and techniques of the functional-analytical approach to
some extremal problems of convex geometry that stem from the Queen Dido problem.
\end{abstract}
\maketitle

\date{October 25, 2023}



\section*{Prologue}
David Mumford,
one of the most beautiful mathematical minds of today, remarked once 
that he honestly carried out some ghastly but wholly straightforward 
calculations while checking something. 
``It took me several hours to do every bit and as I was no wiser at the end... 
I shall omit details here.'' Narrating this episode in \cite{Manin}, another outstanding mathematician, 
Yuri Manin,  concluded: ``The moral: a good proof is one which makes us wiser.''

The following slight abstraction of this thesis transpires: {\it In science
we appraise and appreciate that which makes us wiser.} The notions of
a good theory open up new possibilities of solving particular problems.
Rewarding is the problem whose solution paves way to new fruitful
concepts and methods.  We will shortly present the main ideas and tools
that relate to the classical extremal problems of euclidean geometry.

\section*{The Dido Problem}
The Dido problem is usually acknowleged as the start of  the theory of 
extremal problems \cite{Kelvin}. Dido was a mythical Phoenician Princess. 
Virgil told about the escape of
Dido from her treacherous brother in the first chapter of {\it The Aeneid}.
Dido had to decide about the choice of a~ tract of land near
the future city of Carthage, while satisfying
the famous constraint of selecting ``a space of ground, which (Byrsa call'd,
from the bull's hide) they first inclos'd.''

By the legend, Phoenicians cut the oxhide into thin strips and enclosed
a~large expanse. Now it is customary to think that
the decision by Dido was reduced to
the isoperimetric problem  of finding a figure of   greatest area
among those surrounded by  a curve whose length is given.
It is not excluded that Dido and her subjects solved
the practical versions of the problem when the tower was to be located at the
sea coast and part of the boundary coastline of the tract was somehow prescribed
in advance.


\section*{Minkowski Duality}
Let $\overline{E}:=E\cup\{+\infty\}\cup\{-\infty\}$.
Assume that $H\subset E$ is  a~(convex) cone in $E$,
and so $-\infty$ lies beyond~$H$. A subset $U$  of~ $H$
is {\it convex relative to~}
$H$ or $H$-{\it convex\/}
provided that $U$ is the $H$-{\it support set\/}
$U^H_p:=\{h\in H:h\le p\}$ of some element $p$ of $\overline{E}$.

Alongside the $H$-convex sets we consider
the so-called $H$-convex elements. An element   $p\in \overline{E}$
is  $H$-{\it convex} provided that $p=\sup U^H_p$.
The $H$-convex elements comprise the cone which is denoted by
$\mathscr C(H,\overline{E}$).  We may omit  the references to $H$ when $H$ is clear
from the context. 

Convex elements and sets are ``glued together''
by the {\it Minkowski duality\/} $ \varphi:p\mapsto U^H_p$.
This duality enables us to study convex elements and sets simultaneously \cite{KR}.

\section*{The Space of Convex Bodies}
The Minkowski duality makes $\mathscr V_N$ into a~cone
in the space $C(S_{N-1})$  of continuous functions on the Euclidean unit sphere
$S_{N-1}$, the boundary of the unit ball $\mathfrak z_N$.
This yields
the so-called {\it  Minkowski structure\/} on $\mathscr V_N$.
Addition of the support functions
of convex figures amounts to taking their algebraic sum, also called the
{\it Minkowski addition}. It is worth observing that the
{\it linear span\/}
$[\mathscr V_N]$ of~$\mathscr V_N$ is dense in $C(S_{N-1})$, bears
a~natural structure of a~vector lattice
and is usually referred to as the {\it space of convex sets}.

The study of this space stems from the pioneering breakthrough of
Alexandrov in 1937 and the further insights of
Radstr\"{o}m, H\"{o}rmander, and Pinsker.

\section*{Linear Inequalities over Convex Surfaces}
Reshetnyak demonstrated in his thesis  \cite{Reshetnyak} of 1954 how to
produce inequalities over convex syrfaces.

 A measure $\mu$ {\it linearly majorizes\/} or {\it dominates\/}
a~measure $\nu$  on $S_{N-1}$ provided that to each decomposition of
$S_{N-1}$ into finitely many disjoint Borel sets $U_1,\dots,U_m$
there are measures $\mu_1,\dots,\mu_m$ with sum $\mu$
such that every difference $\mu_k - \nu|_{U_k}$
annihilates all restrictions to $S_{N-1}$ of linear functionals over
$\mathbb R^N$. In symbols, $\mu\,{\gg}{}_{\mathbb R^N} \nu$.

For all sublinear $p$  on  $\mathbb R^N$ we have
$$
\int\limits_{S_N-1} p d\mu \ge  \int\limits_{S_N-1} p d\nu
$$
if   $\mu\,{\gg}{}_{\mathbb R^N} \nu$.

\section*{Choquet Order}
Loomis introduces the concept of affine majorization in \cite{Loomis}:

A~measure $\mu$ {\it affinely majorizes\/} or {\it dominates\/}
a measure $\nu$, both given     on a compact convex subset $Q$ of a locally convex space $X$,
provided that     to each decomposition of
$\nu$ into finitely many summands
$\nu_1,\dots,\nu_m$  there are measures $\mu_1,\dots,\mu_m$
whose sum is $\mu$ and for which every difference
$\mu_k - \nu_k$ annihilates all restrictions
to  $Q$  of affine   functionals over $X$.
In symbols, $\mu\,{\gg}{}_{\Aff(Q)} \nu$.

Cartier, Fell, and Meyer  proved in \cite{CFM} that
$$
\int\limits_{Q} f d\mu \ge  \int\limits_{Q} f d\nu
$$
for each continuous convex function  $f$
on  $Q$   if and only if   $\mu\,{\gg}{}_{\Aff(Q)} \nu$.
An analogous necessity part for linear majorization was published
in \cite{Positive}.

 \section*{Decomposition Theorem}
Majorization is a vast subject. We will use the result in \cite{Kutateladze}:

Assume that $H_1,\dots,H_N$ are cones in a Riesz space~$X$,
while $f$ and $g$ are positive functionals on~$X$.

{\sl
The inequality
$$
f(h_1\vee\dots\vee h_N)\ge g(h_1\vee\dots\vee h_N)
$$
holds for all
$h_k\in H_k$ $(k:=1,\dots,N)$
if and only if to each decomposition
of~$g$ into a~sum of~$N$ positive terms
$g=g_1+\dots+g_N$
there is a decomposition of~$f$ into a~sum of~$N$
positive terms $f=f_1+\dots+f_N$
such that
$$
f_k(h_k)\ge g_k(h_k)\quad
(h_k\in H_k;\ k:=1,\dots,N).
$$
}

\section*{Alexandrov Measures}
  The celebrated {\it Alexandrov Theorem\/} proves the unique existence of
a translate of a convex body given its surface area function; see \cite{AD}.
Each surface area function is an {\it Alexandrov measure}.
So we call a positive measure on the unit sphere which is supported by
no great hypersphere and which annihilates
singletons.

Each Alexandrov measure is a translation-inva\-riant
additive functional over the cone
$\mathscr V_N$.
The cone of positive translation-invariant measures in the
dual $C'(S_{N-1})$ of
 $C(S_{N-1})$ is denoted by~$\mathscr A_N$.

\section*{Blaschke's Sum}
 Given $\mathfrak x, \mathfrak y\in \mathscr V_N$,  the record
$\mathfrak x\,{=}{}_{\mathbb R^N}\mathfrak y$ means that $\mathfrak x$
and $\mathfrak y$  are  equal up to translation or, in other words,
are translates of one another.
So, ${=}{}_{\mathbb R^N}$ is the associate equivalence of
the preorder ${\ge}{}_{\mathbb R^N}$ on $\mathscr V_N$ of
the possibility of inserting one figure into the other
by translation.

The sum of the surface area measures of
$\mathfrak x$ and $\mathfrak y$ generates the unique class
$\mathfrak x\# \mathfrak y$ of translates which is referred to as the
{\it Blaschke sum\/} of $\mathfrak x$ and~$\mathfrak y$; cp.~\cite{Firey}.

\section*{The Natural Duality}
Let $C(S_{N-1})/\mathbb R^N$ stand for the factor space of
$C(S_{N-1})$ by the subspace of all restrictions of linear
functionals on $\mathbb R^N$ to $S_{N-1}$.
Let $[\mathscr A_N]$ be the space $\mathscr A_N-\mathscr A_N$
of translation-invariant measures, in fact, the linear span
of the set of Alexandrov measures.

$C(S_{N-1})/\mathbb R^N$ and $[\mathscr A_N]$ are made dual
by the canonical bilinear form
$$
\gathered
\langle f,\mu\rangle=\frac{1}{N}\int\limits_{S_{N-1}}fd\mu\\
(f\in C(S_{N-1})/\mathbb R^N,\ \mu \in[\mathscr A_N]).
\endgathered
$$

For $\mathfrak x\in\mathscr V_N/\mathbb R^N$ and $\mathfrak y\in\mathscr A_N$,
the quantity
$\langle {\mathfrak x},{\mathfrak y}\rangle$ coincides with the
{\it mixed volume\/}
$V_1 (\mathfrak y,\mathfrak x)$.

\section*{Cones of Feasible Directions}
  Given a cone $K$
in a vector space $X$ in duality with another vector space
$Y$,  the {\it dual\/} of $K$ is
$$
K^*:=\{y\in Y\mid (\forall x\in K)\ \langle x,y\rangle\ge 0\}.
$$

To a convex subset $U$ of $X$ and $\bar x\in U$
there corresponds
$$
U_{\bar x}:=\Fd (U,\bar x):=\{h\in X\mid (\exists \alpha \ge 0)\ \bar x+\alpha h\in U \},
$$
the {\it cone of feasible directions\/}
of $U$ at $\bar x$.

Let $\bar {\mathfrak x}\in{\mathscr A}_N$.
Then the dual  $\mathscr A^*_{N,\bar{\mathfrak x}}$ of the cone of
feasible directions of $\mathscr A_Nn$
at~$\bar{\mathfrak x}$ may be represented as follows
$$
{\mathscr A}^*_{N,\bar{\mathfrak x}}=\{f\in{\mathscr A}^*_N\mid
\langle\bar {\mathfrak x},f\rangle=0\}.
$$

\section*{Dual Cones in Spaces of Surfaces}
Let $\mathfrak x$ and $\mathfrak y$ be  convex figures. Then

$(1)$ $\mu(\mathfrak x)- \mu(\mathfrak y)\in \mathscr V^*_N
\leftrightarrow \mu(\mathfrak x)\,{\gg}{}_{\mathbb R^N} \mu(\mathfrak y)$;

$(2)$ If $\mathfrak x\ge{}_{\mathbb R^N}\mathfrak y$
then  $\mu(\mathfrak x)\,{\gg}{}_{\mathbb R^N} \mu(\mathfrak y)$;

$(3)$ $\mathfrak x\ge{}_{\mathbb R^2}\mathfrak y\leftrightarrow
\mu(\mathfrak x)\,{\gg}{}_{\mathbb R^2} \mu(\mathfrak y)$;

$(4)$ If $\mathfrak y-{\bar{\mathfrak x}}\in\mathscr A^*_{N,\bar{\mathfrak x}}$ then
$\mathfrak y=_{\mathbb R^N}\bar{\mathfrak x}$;

$(5)$ If $\mu (\mathfrak y)-\mu (\bar{\mathfrak x})\in\mathscr V^*_{N,\bar{\mathfrak x}}$
then
$\mathfrak y=_{\mathbb R^N}\bar{\mathfrak x}$.

It stands to reason to avoid  discriminating between a  convex figure,
the respective equivalence class of translates in  $\mathscr V_N/\mathbb R^N$,
and the corresponding measure in $\mathscr A_N$.

\section*{Comparison Between the Structures}

\begin{center}
\includegraphics[scale=.7]{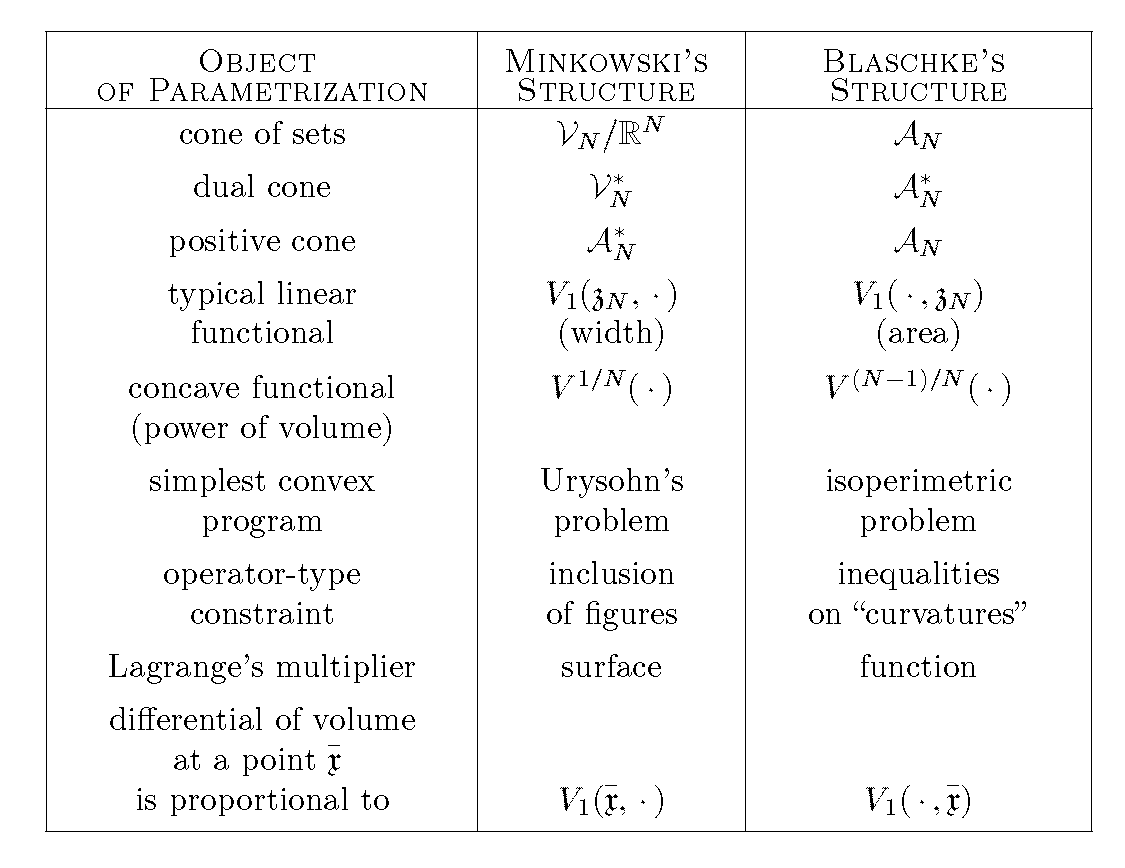}
\end{center}

As regards Brunn--Minkowski theory, see for instance \cite{S} and \cite{G}.

\section*{The External Urysohn Problem}
Among the convex figures, circumscribing $\mathfrak x_0 $ and having
integral breadth fixed, find a convex body of greatest volume.
This is the celebrated Urysohn problem \cite{Urysohn} aggravated
by a  contraint on the location; cp. \cite{R1} and \cite{R2}.

{\sl
A feasible convex body $\bar {\mathfrak x}$ is a solution
to~the external Urysohn problem
if and only if there are a positive  measure~$\mu $
and a positive real $\bar \alpha \in \mathbb R_+$ satisfying

$(1)$ $\bar \alpha \mu
(\mathfrak z_N)\,{\gg}{}_{\mathbb R^N}\mu (\bar {\mathfrak x})+\mu $;

$(2)$~$V(\bar {\mathfrak x})+\frac{1}{N}\int_{S_{N-1}}
\bar {\mathfrak x}d\mu =\bar \alpha V_1 (\mathfrak z_N,\bar {\mathfrak x})$;

$(3)$~$\bar{\mathfrak x}(z)={\mathfrak x}_0 (z)$
for all $z$ in the support of~$\mu $.
}

\section*{Solutions}
  If ${\mathfrak x}_0 ={\mathfrak z}_{N-1}$ then $\bar{\mathfrak x}$
is a {\it spherical lens} and $\mu$ is the restriction
of the surface area function
of the ball of radius
$\bar \alpha ^{1/(N-1)}$
to the complement of the support of the lens to~$S_{N-1}$.

If ${\mathfrak x}_0$ is an equilateral triangle then the solution
$\bar {\mathfrak x}$ looks as follows:

\centerline{\includegraphics[scale=1.0]{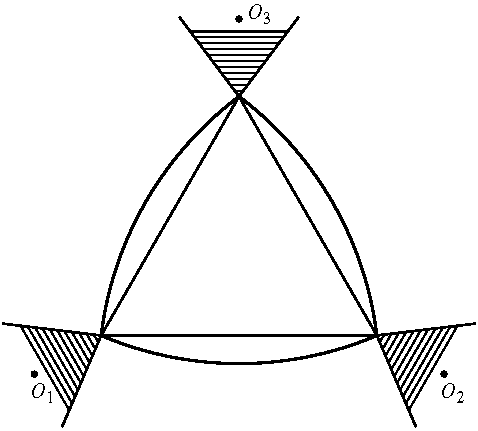}}

$\bar {\mathfrak x}$ is the union of~${\mathfrak x}_0$
and  three congruent slices of a circle of radius~$\bar \alpha$ and
centers $O_1$--$O_3$, while
$\mu$ is the restriction of $\mu(\mathfrak z_2)$
to the subset of $S_1$ comprising the endpoints
of the unit vectors of the shaded zone.

\section*{Symmetric Solutions}
This is  the general solution of the internal Urysohn problem inside a triangle
in the class of centrally symmetric convex figures:

\centerline{\includegraphics[scale=1.0]{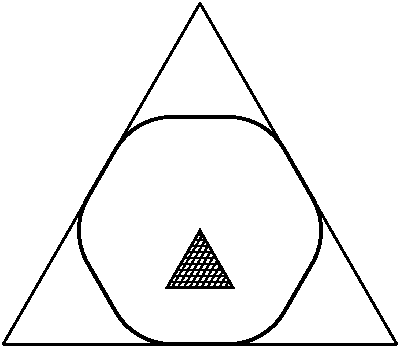}}

\section*{Current Hyperplanes}
  Find two convex figures $\bar{\mathfrak x}$ and $\bar{\mathfrak y}$
lying in a given convex body
$\mathfrak x_o$,
 separated by a~hyperplane with the unit outer normal~$z_0$,
and having the greatest total volume
of $\bar{\mathfrak x}$ and~$\bar{\mathfrak y}$
given the sum of their integral breadths.

{\sl
A feasible pair of convex bodies $\bar{\mathfrak x}$ and $\bar{\mathfrak y}$
solves the internal Urysohn problem with a current hyperplane
if and only if
there are convex figures  $\mathfrak x$ and $\mathfrak y$
and positive reals
$\bar\alpha $ and $\bar\beta$  satisfying

{\rm(1)}  $\bar{\mathfrak x}=\mathfrak x \# \bar\alpha\mathfrak z_N$;

{\rm(2)} $\bar{\mathfrak y}=\mathfrak y \# \bar\alpha\mathfrak z_N$;

{\rm(3)} $\mu(\mathfrak x) \ge \bar\beta\varepsilon_{z_0} $, $\mu(\mathfrak y) \ge \bar\beta\varepsilon_{-z_0} $;

{\rm(4)} $\bar {\mathfrak x}(z)=\mathfrak x_0 (z)$ for all $z\in \spt(\mathfrak x)\setminus \{z_0\} $;

{\rm(5)} $\bar {\mathfrak y}(z)=\mathfrak x_0 (z)$ for all $z\in \spt(\mathfrak x)\setminus \{-z_0\} $,
\noindent
with $\spt(\mathfrak x)$ standing for the {\it support\/} of $\mathfrak x$,
i.e. the support of the surface area measure $\mu(\mathfrak x)$
of~$\mathfrak x$.
}

\section*{ Pareto Optimality}
  Consider a~bunch of  economic agents
each of which intends to maximize his own income.
The {\it Pareto efficiency principle\/}  asserts
that  as an effective agreement of the conflicting goals it is reasonable
to take any state in which nobody can increase his income in any way other
than diminishing the income of at least one of the other fellow members.

Formally speaking, this implies the search of the maximal elements
of the set comprising the tuples of incomes of the agents
at every state; i.e., some vectors of a finite-dimensional
arithmetic space endowed with the coordinatewise order. Clearly,
the concept of Pareto optimality was
already abstracted to arbitrary
ordered vector spaces; see, for instance, \cite{Z}--\cite{Subdif}.

\section*{ Vector Isoperimetric Problem}
 Given are some convex bodies
$\mathfrak y_1,\dots,\mathfrak y_M$.
Find a convex body $\mathfrak x$ encompassing a given volume
and minimizing each of the mixed volumes $V_1(\mathfrak x,\mathfrak y_1),\dots,V_1(\mathfrak x,\mathfrak y_M)$.
In symbols,
$$
\mathfrak x\in\mathscr A_N;\
\widehat p(\mathfrak x)\ge \widehat p(\bar{\mathfrak x});\
(\langle\mathfrak y_1,\mathfrak x\rangle,\dots,\langle\mathfrak y_M,\mathfrak x\rangle)\rightarrow\inf\!.
$$
Clearly, this is a~Slater regular convex program in the Blaschke structure.

 {\sl
Each Pareto-optimal solution $\bar{\mathfrak x}$ of the vector isoperimetric problem
has the form}
$$
\bar{\mathfrak x}=\alpha_1{\mathfrak y}_1+\dots+\alpha_m{\mathfrak y}_m,
$$
where $\alpha_1,\dots,\alpha_m$ are positive reals.

\section*{The Leidenfrost Problem}
Given the volume of a three-dimensional
convex figure, minimize its  surface area and vertical  breadth.

By symmetry everything reduces to an analogous plane two-objective problem,
whose every Pareto-optimal solution is by~2 a~{\it stadium\/},
a weighted Minkowski sum of a disk and
a horizontal straight line segment.

{\sl A plane spheroid, a Pareto-optimal solution of the Leidenfrost problem,
is the result of rotation of a stadium around the vertical axis through
the center of the stadium}.

\section*{ Internal Urysohn Problem with Flattening}
Given are some~convex body
$\mathfrak x_0\in\mathscr V_N$ and some flattening direction~ $\bar z\in S_{N-1}$.
Considering $\mathfrak x\subset\mathfrak x_0$ of
fixed integral breadth, maximize the volume of~$\mathfrak x$ and  minimize the
breadth of $\mathfrak x$ in the flattening direction:
$\mathfrak x\in\mathscr V_N;\
\mathfrak x\subset{\mathfrak x}_0;\
\langle \mathfrak x,{\mathfrak z}_N\rangle \ge \langle\bar{\mathfrak x},{\mathfrak z}_N\rangle;\
(-p(\mathfrak x), b_{\bar z}(\mathfrak x)) \to\inf\!.
$

{\sl For a feasible convex body $\bar{\mathfrak x}$ to be Pareto-optimal in
the internal Urysohn problem with the flattening
direction~$\bar z$ it is necessary and sufficient that there be
positive reals $\alpha, \beta$ and a~convex figure $\mathfrak x$ satisfying}
$$
\gathered
\mu(\bar{\mathfrak x})=\mu(\mathfrak x)+ \alpha\mu({\mathfrak z}_N)+\beta(\varepsilon_{\bar z}+\varepsilon_{-\bar z});\\
\bar{\mathfrak x}(z)={\mathfrak x}_0(z)\quad (z\in\spt(\mu(\mathfrak x)).
\endgathered
$$

\section*{Rotational Symmetry}
Assume that a plane convex figure ${\mathfrak x}_0\in\mathscr V_2$ has the symmetry axis $A_{\bar z}$
with generator~$\bar z$.  Assume further that ${\mathfrak x}_{00}$ is the result of rotating
$\mathfrak x_0$  around the symmetry axis $A_{\bar z}$ in~$\mathbb R^3$.
$$
\gathered
\mathfrak x\in\mathscr V_3;\\
\mathfrak x  \text{\ is\ a\ convex\ body\ of\ rotation\ around}\ A_{\bar z};\\
\mathfrak x\supset{\mathfrak x}_{00};\
\langle {\mathfrak z}_N, \mathfrak x\rangle \ge \langle{\mathfrak z}_N,\bar{\mathfrak x}\rangle;\\
(-p(\mathfrak x), b_{\bar z}(\mathfrak x)) \to\inf\!.
\endgathered
$$

{\sl Each Pareto-optimal solution  is the result
of rotating around the symmetry axis a Pareto-optimal solution of the plane internal
Urysohn problem with flattening in the direction of the axis}.

\section*{Soap Bubbles}
Little is known about the analogous problems in arbitrary dimensions; 
cp. \cite{I} and \cite{L}.
An especial place
is occupied by the result of Porogelov   who
demonstrated that the ``soap bubble'' in a tetrahedron
has the form of the result of the rolling of a ball over a~solution
of the internal Urysohn problem, i.e. the weighted Blaschke sum of
a tetrahedron and a ball; see \cite{Pogorelov}.

\section*{Is Dido's Problem Solved?}
  From a utilitarian standpoint, the answer is
definitely in the affirmative. There is no evidence that Dido
experienced any difficulties, showed indecisiveness, and procrastinated the choice of the tract of land.
Practically speaking, the situation in which Dido  made her decision
was not as primitive as it seems at the first glance.
The appropriate generality is unavailable in the
mathematical model known as the classical isoperimetric problem.
For more details see \cite{Multi1} and \cite{Multi2}.

Dido's problem inspiring our ancestors remains the same intellectual
challenge as Kant's starry heavens above and moral law within.

\bibliographystyle{plain}

\end{document}